\documentclass[12pt]{article}
\usepackage{amssymb}

\newcommand{\be}{\begin{equation}}
\newcommand{\ee}{\end{equation}}
\newcommand{\bea}{\begin{eqnarray}}
\newcommand{\eea}{\end{eqnarray}}
\newcommand{\barray}{\begin{array}}
\newcommand{\earray}{\end{array}}
\newcommand{\pa}{\partial}
\newcommand{\nn}{\nonumber}
\newcommand{\bitem}{\begin{itemize}}
\newcommand{\eitem}{\end{itemize}}
\newtheorem{teo}{Theorem}[section]
\newcommand{\bt}{\begin{teo}}
\newcommand{\et}{\end{teo}}
\newtheorem{Def}{Definition}[section]
\newcommand{\bd}{\begin{Def}}
\newcommand{\ed}{\end{Def}}
\newtheorem{lem}{Lemma}[section]
\newcommand{\bl}{\begin{lem}}
\newcommand{\el}{\end{lem}}
\newtheorem{prop}{Proposition}[section]
\newcommand{\bp}{\begin{prop}}
\newcommand{\ep}{\end{prop}}
\newtheorem{cor}{Corollary}[section]
\newcommand{\bc}{\begin{cor}}
\newcommand{\ec}{\end{cor}}
\newtheorem{ex}{Example}[section]
\newcommand{\bex}{\begin{ex}}
\newcommand{\eex}{\end{ex}}
\newtheorem{rem}{Remark}[section]
\newcommand{\br}{\begin{rem}}
\newcommand{\er}{\end{rem}}

\catcode `\@=11
\@addtoreset{equation}{section}

\begin{document}

{53B25, 53D45, 35Q58, 81T45}
\bigskip

\begin{center}
{\Large \textbf{Theory of Submanifolds,
Associativity Equations in 2D Topological Quantum Field Theories,
and Frobenius Manifolds\footnote{This research was supported by the
Max-Planck-Institut f\"{u}r Mathematik (Bonn, Germany),
the Russian Science Support Foundation,
the Russian Foundation for Basic Research
(Grant No. 05-01-00170), and the Program of Support for
Leading Scientific Schools (Grant No. NSh-4182.2006.1).}}}
\end{center}

\medskip

\begin{center}
{\large {O. I. Mokhov}}
\end{center}
\begin{center}
{Center for Nonlinear Studies,
L.D.Landau Institute for Theoretical Physics,\\
Russian Academy of Sciences,
Kosygina 2, Moscow, 117940, Russia\\
Department of Geometry and Topology,
Faculty of Mechanics and Mathematics, \\
M.V.Lomonosov Moscow State University,
Moscow, 119992, Russia\\
e-mail: mokhov@mi.ras.ru; mokhov@landau.ac.ru; mokhov@bk.ru}
\end{center}

\begin{flushright}
{To memory of my wonderful mother \hspace {15mm} \ \\
Maya Nikolaevna Mokhova
(04.05.1926 -- 12.09.2006)}
\end{flushright}

\begin{abstract}
We prove that the associativity
equations of two-dimensional topological quantum field theories
are very natural reductions of the fundamental nonlinear equations
of the theory of submanifolds in
pseudo-Euclidean spaces and give a natural class of {\it potential}
flat torsionless submanifolds. We show that
all potential flat torsionless submanifolds
in pseudo-Euclidean spaces bear natural structures of Frobenius algebras
on their tangent spaces. These Frobenius structures are
generated by the corresponding flat first fundamental form
and the set of the second fundamental forms of the submanifolds
(in fact, the structural constants are given
by the set of the Weingarten operators of the submanifolds).
We prove in this paper
that
each $N$-dimensional Frobenius manifold can locally be represented as a
potential flat torsionless submanifold in a $2N$-dimensional pseudo-Euclidean
space. By our construction this submanifold is uniquely determined
up to motions.
Moreover, in this paper we consider a nonlinear system, which is a
natural generalization of the associativity equations,
namely, the system describing all flat torsionless submanifolds
in pseudo-Euclidean spaces, and prove that this system is integrable by
the inverse scattering method.
\end{abstract}

\section{Introduction. Associativity equations and \\
Frobenius structures} \label{introduction}

In this paper we prove that the associativity
equations of two-dimensional topological quantum field theories
(the Witten--Dijkgraaf--Verlinde--Verlinde and Dubrovin equations,
see \cite{1}) for
a function ({\it potential} or {\it prepotential})
$\Phi = \Phi (u^1, \ldots, u^N)$,
\be
\sum_{k = 1}^N \sum_{l = 1}^N
{\pa^3 \Phi \over \pa u^i \pa u^j \pa u^k} \eta^{kl}
{\pa^3 \Phi \over \pa u^l \pa u^m \pa u^n} =
\sum_{k = 1}^N \sum_{l = 1}^N
{\pa^3 \Phi \over \pa u^i \pa u^m \pa u^k} \eta^{kl}
{\pa^3 \Phi \over \pa u^l \pa u^j \pa u^n},  \label{ass1}
\ee
where
$\eta^{ij}$ is an arbitrary constant nondegenerate symmetric matrix,
$\eta^{ij} = \eta^{ji},$ $\eta^{ij} = {\rm const},$
$\det (\eta^{ij}) \neq 0$,
are very natural reductions of the fundamental nonlinear equations
of the theory of submanifolds in
pseudo-Euclidean spaces and give a natural class of {\it potential}
flat torsionless submanifolds. All potential flat torsionless submanifolds
in pseudo-Euclidean spaces bear natural structures of Frobenius algebras
on their tangent spaces. These Frobenius structures are
generated by the corresponding flat first fundamental form
and the set of the second fundamental forms of the submanifolds
(in fact, the structural constants are given
by the set of the Weingarten operators of the submanifolds).
We recall that each solution
$\Phi (u^1, \ldots, u^N)$ of the
associativity equations (\ref{ass1}) gives $N$-parametric deformations of
Frobenius algebras, i.e., commutative associative algebras equipped by
nondegenerate invariant symmetric bilinear forms.
Indeed, consider algebras $A (u)$ in an $N$-dimensional vector space with
basis $e_1, \ldots, e_N$ and multiplication (see \cite{1})
\be
e_i \circ e_j = c^k_{ij} (u) e_k, \ \ \ \
c^k_{ij} (u) = \eta^{ks} {\pa^3 \Phi \over \pa u^s \pa u^i \pa u^j} e_k.
\label{al1}
\ee
For all values of the parameters $u = (u^1, \ldots, u^N)$ the algebras
$A (u)$ are commutative, $e_i \circ e_j =
e_j \circ e_i,$ and the associativity condition
\be
(e_i \circ e_j) \circ e_k = e_i \circ (e_j \circ e_k)  \label{al2}
\ee
in the algebras $A (u)$ is equivalent to equations (\ref{ass1}).
The matrix $\eta_{ij}$ inverse to the matrix $\eta^{ij}$,
$\eta^{is} \eta_{sj} = \delta^i_j$, defines a
nondegenerate invariant symmetric bilinear form on the algebras
$A (u)$,
\be
\langle e_i, e_j \rangle  = \eta_{ij}, \ \ \ \
\langle e_i \circ e_j, e_k \rangle =
\langle e_i, e_j \circ e_k \rangle.  \label{al3}
\ee
Recall that locally the tangent space at every point of
a Frobenius manifold (see \cite{1}) bears the Frobenius algebra structure
(\ref{al1})--(\ref{al3}), which is defined by a solution of the
associativity equations (\ref{ass1}) and smoothly depends on the point.
One should also impose additional conditions on Frobenius manifolds,
but we do not consider these conditions here. We prove in this paper
that
each $N$-dimensional Frobenius manifold can locally be represented as a
potential flat torsionless submanifold in a $2N$-dimensional pseudo-Euclidean
space. By our construction this submanifold is uniquely determined
up to motions.
Moreover, in this paper we consider a nonlinear system, which is a
natural generalization of the associativity equations (\ref{ass1}),
namely, the system describing all flat torsionless submanifolds
in pseudo-Euclidean spaces, and prove that this system is integrable by
the inverse scattering method. The connection of the construction
with integrable hierarchies, nonlocal Hamiltonian operators of hydrodynamic
type with
flat metrics, Poisson pencils and recursion operators can be found in
\cite{2}. The results of application of this construction to the theory
of Frobenius manifolds will be published in a separate paper.

\section{Fundamental nonlinear equations of the theory of submanifolds in
Euclidean spaces} \label{section1}

Let us consider an arbitrary smooth $N$-dimensional submanifold
$M^N$ in an $(N +~L)$-dimensional Euclidean space $E^{N + L}$,
$M^N \subset E^{N + L}$,
and introduce the standard classic notation. Let the submanifold
$M^N$ be given locally by a smooth vector-function
$ r (u^1, \ldots, u^N)$ of $N$ independent variables
$(u^1, \ldots, u^N)$ (some independent parameters on the submanifold),
$ r (u^1, \ldots, u^N) =
(z^1 (u^1, \ldots, u^N), \ldots, z^{N + L} (u^1, \ldots, u^N)),$
where $(z^1, \ldots, z^{N + L})$ are coordinates in the
Euclidean space $E^{N + L}$, $(z^1, \ldots, z^{N + L}) \in E^{N + L}$,
$(u^1, \ldots, u^N)$ are local coordinates (parameters) on $M^N$,
${\rm rank} (\pa z^i/\pa u^j) = N$ (here $1 \leq i \leq N + L,$
$1 \leq j \leq N$).
Then ${\pa r/ \pa u^i} = r_{u^i},$ $1 \leq i \leq N,$ are
tangent vectors at an arbitrary point $u = (u^1, \ldots, u^N)$ on $M^N$.
Let ${\bf N}_u$ be the
normal space of the submanifold
$M^N$ at an arbitrary point
$u = (u^1, \ldots, u^N)$ on $M^N$,
$N_u = \langle n_1, \ldots, n_L \rangle$, where
$n_{\alpha}$, $1 \leq \alpha \leq L,$ are orthonormalized normals,
$(n_{\alpha}, r_{u^i}) = 0,$ $1 \leq \alpha \leq L,$
$1 \leq i \leq N,$
$(n_{\alpha}, n_{\beta}) = 0,$
$1 \leq \alpha, \beta \leq L,$ $\alpha \neq \beta,$ and
$(n_{\alpha}, n_{\alpha}) = 1,$ $1 \leq \alpha \leq L.$

Then ${\bf I} = d s^2 = g_{ij} (u) d u^i d u^j,$
$g_{ij} (u) = (r_{u^i}, r_{u^j}),$
is the first fundamental form and
${\bf II}_{\alpha} = \omega_{\alpha, ij} (u) d u^i d u^j,$
$\omega_{\alpha, ij} (u) = (n_{\alpha}, r_{u^i u^j}),$
$1 \leq \alpha \leq L,$ are the second fundamental forms
of the submanifold $M^N$.

Since
the set of vectors
$(r_{u^1} (u), \ldots, r_{u^N} (u), n_1 (u), \ldots, n_L (u))$
forms a basis in $E^{N + L}$ at each point of the submanifold $M^N$, we
can decompose each of the vectors $n_{\alpha, u^i} (u),$
$1 \leq \alpha \leq L,$
$1 \leq i \leq N,$ with respect to this basis, namely,
$n_{\alpha, u^i} (u) = A^k_{\alpha, i} (u) r_{u^k} (u) +
\varkappa_{\alpha \beta, i} (u) n_{\beta} (u)$, where $A^k_{\alpha, i} (u)$
and $\varkappa_{\alpha \beta, i} (u)$ are some coefficients
depending on $u$ (the {\it Weingarten decomposition}).
It is easy to prove that
$A^k_{\alpha, i} (u) = - \omega_{\alpha, ij} (u) g^{jk} (u),$
where $g^{jk} (u)$ is the contravariant metric inverse to the
first fundamental form $g_{ij} (u)$, $g^{is} (u) g_{sj} (u) = \delta^i_j.$
The coefficients $\varkappa_{\alpha \beta, i} (u)$ are said to be
the {\it coefficients of torsion of the submanifold $M^N$,}
$\varkappa_{\alpha \beta, i} (u) = (n_{\alpha, u^i} (u), n_{\beta} (u)).$
It is also easy to prove that the coefficients
$\varkappa_{\alpha \beta, i} (u)$ are skew-symmetric with respect to
the indices $\alpha$ and $\beta$,
$\varkappa_{\alpha \beta, i} (u) = - \varkappa_{\beta \alpha, i} (u)$,
and form covariant tensors (1-forms) with respect to the index $i$ on
the submanifold $M^N$.
The 1-forms $\varkappa_{\alpha \beta, i} (u) d u^i$ are said to be
the {\it torsion forms of the submanifold $M^N$}.

It is well known that for each submanifold $M^N$ the forms $g_{ij} (u)$,
$\omega_{\alpha, ij} (u)$ and $\varkappa_{\alpha \beta, i} (u)$
satisfy the Gauss equations, the Codazzi equations and the Ricci equations,
which are the fundamental equations of the theory of submanifolds.
The Gauss equations have the form
\be
R_{ijkl} (u) = \sum_{\alpha = 1}^L
\left (\omega_{\alpha, jl} (u) \omega_{\alpha, ik} (u) -
\omega_{\alpha, jk} (u) \omega_{\alpha, il} (u) \right ), \label{g1}
\ee
where $R_{ijkl} (u)$ is the tensor of Riemannian curvature of
the first fundamental form $g_{ij} (u)$. The Codazzi equations
have the form
\be
\nabla_k (\omega_{\alpha, ij} (u)) -
\nabla_j (\omega_{\alpha, ik} (u)) = \varkappa_{\alpha \beta, k} (u)
\omega_{\beta, ij} (u) - \varkappa_{\alpha \beta, j} (u)
\omega_{\beta, ik} (u), \label{c1}
\ee
where $\nabla_k$ is the covariant differentiation generated
by the Levi-Civita connection of the first fundamental form $g_{ij} (u)$.
The Ricci equations have the form
\bea
&&
\nabla_k (\varkappa_{\alpha \beta, i} (u)) -
\nabla_i (\varkappa_{\alpha \beta, k} (u)) +
\sum_{\gamma = 1}^L \left (\varkappa_{\alpha \gamma, i} (u)
\varkappa_{\gamma \beta, k} (u) - \varkappa_{\alpha \gamma, k} (u)
\varkappa_{\gamma \beta, i} (u) \right ) + \nn\\
&&
+ \left (\omega_{\alpha, kl} (u) \omega_{\beta, ji} (u) - \omega_{\alpha, il} (u)
\omega_{\beta, jk} (u) \right ) \, g^{lj} (u) = 0. \label{r1}
\eea

{\bf The Bonnet theorem}.
{\it Let $K^N$ be an arbitrary smooth $N$-dimensional Riemannian
manifold with a metric
$g_{ij} (u) du^i du^j$. Let some 2-forms
$\omega_{\alpha, ij} (u) du^i du^j$, $1 \leq \alpha \leq L,$
and some 1-forms $\varkappa_{\alpha \beta, i} (u)$,
$1 \leq \alpha, \beta \leq L,$ be given in a simply connected
domain of the manifold $K^N$. If $\omega_{\alpha, ij} (u) =
\omega_{\alpha, ji} (u)$, $\varkappa_{\alpha \beta, i} (u) =
- \varkappa_{\beta \alpha, i} (u)$, and the Gauss
equations {\rm (\ref{g1})},
the Codazzi equations {\rm (\ref{c1})}, and the Ricci equations
{\rm (\ref{r1})} are satisfied for
the forms $g_{ij} (u)$,
$\omega_{\alpha, ij} (u)$ and $\varkappa_{\alpha \beta, i} (u)$,
then there exists a unique (up to motions) smooth
$N$-dimensional submanifold $M^N$
in an $(N + L)$-dimensional Euclidean space
$E^{N + L}$ with the first fundamental form
$d s^2 = g_{ij} (u) du^i du^j$, the second fundamental forms
$\omega_{\alpha, ij} (u) d u^i d u^j$ and the torsion forms
$\varkappa_{\alpha \beta, i} (u) d u^i$.}

Similar fundamental equations and the Bonnet theorem are true for all
{\it totally nonisotropic submanifolds in pseudo-Euclidean spaces}
(we recall that if we have a submanifold in an
arbitrary pseudo-Euclidean space $E^m_n$,
then the metric induced on the submanifold from the ambient
pseudo-Euclidean space $E^m_n$ is nondegenerate
if and only if this submanifold is totally nonisotropic).

\section{Flat submanifolds with zero torsion \\
in pseudo-Euclidean spaces} \label{section2}

Let us consider totally nonisotropic smooth  $N$-dimensional
submanifolds with {\it zero torsion} in an $(N + L)$-dimensional
pseudo-Euclidean space, i.e.,
the torsion forms of submanifolds of this class vanish,
$\varkappa_{\alpha \beta, i} (u) = 0$.
In the normal spaces $N_u$
we will also use bases $n_{\alpha}$, $1 \leq \alpha \leq L,$
with arbitrary admissible Gram matrices $\mu_{\alpha \beta}$,
$(n_{\alpha}, n_{\beta}) = \mu_{\alpha \beta}$,
$\mu_{\alpha \beta}= {\rm const}$,
$\mu_{\alpha \beta} = \mu_{\beta \alpha}$, $\det \mu_{\alpha \beta}
\neq 0$ (the signature of the metric $\mu_{\alpha \beta}$ is
determined by the
signature of the first fundamental form of the submanifold
and the signature of the ambient
pseudo-Euclidean space).

For torsionless $N$-dimensional submanifolds in an $(N + L)$-dimensional
pseudo-Euclidean space we have the following system of fundamental equations,
the Gauss equations
\be
R_{ijkl} (u) = \sum_{\alpha = 1}^L \sum_{\beta = 1}^L
\mu^{\alpha \beta}
(\omega_{\alpha, ik} (u) \omega_{\beta, jl} (u) -
\omega_{\alpha, il} (u) \omega_{\beta, jk} (u)),  \label{g2}
\ee
where $\mu^{\alpha \beta}$ is the inverse to the matrix
$\mu_{\alpha \beta}$, $\mu^{\alpha \gamma} \mu_{\gamma \beta} =
\delta^{\alpha}_{\beta}$, the Codazzi equations
\be
\nabla_k (\omega_{\alpha, ij} (u)) =
\nabla_j (\omega_{\alpha, ik} (u)),   \label{c2}
\ee
and the Ricci equations
\be
g^{ij} (u) \, (\omega_{\alpha, ik} (u) \omega_{\beta, jl} (u)
- \omega_{\alpha, il} (u)
\omega_{\beta, jk} (u)) = 0.  \label{r2}
\ee

Now let $g_{ij} (u)$ be a flat metric, i.e., we consider flat
torsionless $N$-dimensional submanifolds $M^N$
in an $(N + L)$-dimensional
pseudo-Euclidean space. Then we can consider that
$u = (u^1, \ldots, u^N)$ are certain flat coordinates of the
metric $g_{ij} (u)$  on $M^N$. In flat coordinates the metric is
a constant nondegenerate symmetric matrix $\eta_{ij}$,
$\eta_{ij} = \eta_{ji},$ $\eta_{ij} = {\rm const},$
$\det (\eta_{ij}) \neq 0$, and
the Codazzi equations (\ref{c2}) have the form
\be
{\pa \omega_{\alpha, ij} \over \pa u^k} =
{\pa \omega_{\alpha, ik} \over \pa u^j}.
\ee
Thus there exist locally some functions $\chi_{\alpha, i} (u),$
$1 \leq \alpha \leq L,$ $1 \leq i \leq N,$ such that
\be
\omega_{\alpha, ij} (u) = {\pa \chi_{\alpha, i} \over \pa u^j}.
\ee
From symmetry of the second fundamental forms
$\omega_{\alpha, ij} (u) = \omega_{\alpha, ji} (u)$
we have
\be
{\pa \chi_{\alpha, i} \over \pa u^j} =
{\pa \chi_{\alpha, j} \over \pa u^i}.
\ee
Therefore, there exist locally some functions $\psi_{\alpha} (u),$
$1 \leq \alpha \leq L,$ such that
\be
\chi_{\alpha, i} (u) = {\pa \psi_{\alpha} \over \pa u^i}, \ \ \ \
\omega_{\alpha, ij} (u) = {\pa^2 \psi_{\alpha} \over \pa u^i \pa u^j}.
\ee
Thus we have proved the following important lemma.
\bl
All the second fundamental forms of each
flat torsionless submanifold
in a pseudo-Euclidean space are Hessians in any flat coordinates
in any simply connected domain on the
submanifold.
\el
Moreover, in any flat coordinates the Gauss equations (\ref{g2})
have the form
\be
\sum_{\alpha = 1}^L \sum_{\beta = 1}^L
\mu^{\alpha \beta}
\left ({\pa^2 \psi_{\alpha} \over \pa u^i \pa u^k}
{\pa^2 \psi_{\beta} \over \pa u^j \pa u^l} -
{\pa^2 \psi_{\alpha} \over \pa u^i \pa u^l}
{\pa^2 \psi_{\beta} \over \pa u^j \pa u^k} \right ) = 0 \label{g3}
\ee
and the Ricci equations (\ref{r2}) have the form
\be
\sum_{i = 1}^N \sum_{j = 1}^N
\eta^{ij}
\left ({\pa^2 \psi_{\alpha} \over \pa u^i \pa u^k}
{\pa^2 \psi_{\beta} \over \pa u^j \pa u^l} -
{\pa^2 \psi_{\alpha} \over \pa u^i \pa u^l}
{\pa^2 \psi_{\beta} \over \pa u^j \pa u^k} \right ) = 0,  \label{r3}
\ee
where $\eta^{ij}$ is the inverse to the matrix
$\eta_{ij}$, $\eta^{is} \eta_{sj} =
\delta^i_j$.

\bt
The class of $N$-dimensional flat torsionless submanifolds
in $(N + L)$-dimensional
pseudo-Euclidean spaces is described (in flat coordinates) by the
system of nonlinear equations {\rm (\ref{g3}), (\ref{r3})} for
functions $\psi_{\alpha} (u),$ $1 \leq \alpha \leq L.$
Here $\eta^{ij}$ and $\mu^{\alpha \beta}$
are arbitrary constant nondegenerate symmetric matrices,
$\eta^{ij} = \eta^{ji},$ $\eta^{ij} = {\rm const},$
$\det (\eta^{ij}) \neq 0$,
$\mu^{\alpha \beta}= {\rm const}$,
$\mu^{\alpha \beta} = \mu^{\beta \alpha}$, $\det \mu^{\alpha \beta}
\neq 0$, the signature of the ambient $(N + L)$-dimensional
pseudo-Euclidean space is the sum of the signatures of
the metrics $\eta^{ij}$ and $\mu^{\alpha \beta}$,
${\bf I} = d s^2 = \eta_{ij} d u^i d u^j$ is the
first fundamental form, where $\eta_{ij}$ is the inverse to the matrix
$\eta^{ij}$, $\eta^{is} \eta_{sj} =
\delta^i_j,$
${\bf II}_{\alpha} =
(\pa^2 \psi_{\alpha} / (\pa u^i \pa u^j)) d u^i d u^j,$
$1 \leq \alpha \leq L,$ are the second fundamental forms
given by the Hessians of the
functions  $\psi_{\alpha} (u),$ $1 \leq \alpha \leq L$.
\et

According to the Bonnet theorem
any solution $\psi_{\alpha} (u),$ $1 \leq \alpha \leq L,$
of the nonlinear system (\ref{g3}), (\ref{r3})
defines a unique (up to motions) $N$-dimensional flat torsionless
submanifold of
the corresponding $(N + L)$-dimensional pseudo-Euclidean space with the
first fundamental form $\eta_{ij} d u^i d u^j$ and the
second fundamental forms
$\omega_{\alpha} (u) =
(\pa^2 \psi_{\alpha} / (\pa u^i \pa u^j))
d u^i d u^j$, $1 \leq \alpha \leq L$, given by the Hessians of the
functions  $\psi_{\alpha} (u),$ $1 \leq \alpha \leq L$.
It is obvious that we can always add arbitrary terms linear in
the coordinates $(u^1, \ldots, u^N)$ to any solution of the system
(\ref{g3}), (\ref{r3}), but the set of the second fundamental forms
and the corresponding submanifold
will be the same. Moreover, any two sets of
the second fundamental forms $\omega_{\alpha, ij} (u) =
\pa^2 \psi_{\alpha} / (\pa u^i \pa u^j)$, $1 \leq \alpha \leq L$,
coincide if and only if the corresponding
functions $\psi_{\alpha} (u),$ $1 \leq \alpha \leq L$,
coincide up to linear terms, so
we must not distinguish solutions of
the nonlinear system (\ref{g3}), (\ref{r3}) up to
terms linear in the coordinates $(u^1, \ldots, u^N)$.

\bt
The nonlinear system {\rm (\ref{g3}), (\ref{r3})} is
integrable by the inverse scattering method.
\et

Consider the following linear problem for
vector-functions $\pa a(u) / \pa u^i$
and $b_{\alpha} (u)$, $1 \leq \alpha \leq L$:
\be
{\pa^2 a \over \pa u^i \pa u^j} = \lambda \, \mu^{\alpha \beta}
\omega_{\alpha, ij} (u)  b_{\beta} (u),\ \ \ \
{\pa b_{\alpha} \over \pa u^i} =
\rho \, \eta^{kj} \omega_{\alpha, ij} (u) {\pa a \over \pa u^k}, \label{t1}
\ee
where $\eta^{ij}$, $ 1 \leq i, j \leq N,$ and $\mu^{\alpha \beta},$
$1 \leq \alpha, \beta \leq L$,
are arbitrary constant nondegenerate symmetric matrices,
$\eta^{ij} = \eta^{ji},$ $\eta^{ij} = {\rm const},$
$\det (\eta^{ij}) \neq 0$,
$\mu^{\alpha \beta}= {\rm const}$,
$\mu^{\alpha \beta} = \mu^{\beta \alpha}$, $\det \mu^{\alpha \beta}
\neq 0$; $\omega_{\alpha, ij} (u)$, $1 \leq \alpha \leq L,$ must be
symmetric matrix functions, $\omega_{\alpha, ij} (u) =
\omega_{\alpha, ji} (u);$
 $\lambda$ and $\rho$ are arbitrary constants (parameters).

The consistency condition for the linear system (\ref{t1})
gives the nonlinear system (\ref{g3}), (\ref{r3})
describing the class of $N$-dimensional flat torsionless submanifolds
in $(N + L)$-dimensional pseudo-Euclidean spaces. Indeed,
we have
\bea
&&
{\pa^3 a \over \pa u^i \pa u^j \pa u^k} = \lambda \, \mu^{\alpha \beta}
{\pa \omega_{\alpha, ij} \over \pa u^k}  b_{\beta} (u) +
\lambda \, \mu^{\alpha \beta}
\omega_{\alpha, ij} (u)  {\pa b_{\beta} \over \pa u^k} = \nn\\
&&
= \lambda \, \mu^{\alpha \beta}
{\pa \omega_{\alpha, ij} \over \pa u^k}  b_{\beta} (u) +
\lambda \, \mu^{\alpha \beta}
\omega_{\alpha, ij} (u)
\rho \, \eta^{ls} \omega_{\beta, ks} (u) {\pa a \over \pa u^l} = \nn\\
&&
= \lambda \, \mu^{\alpha \beta}
{\pa \omega_{\alpha, ik} \over \pa u^j}  b_{\beta} (u) +
\lambda \, \mu^{\alpha \beta}
\omega_{\alpha, ik} (u)
\rho \, \eta^{ls} \omega_{\beta, js} (u) {\pa a \over \pa u^l},
\eea
whence we obtain
\be
{\pa \omega_{\alpha, ij} (u) \over \pa u^k} =
{\pa \omega_{\alpha, ik} (u) \over \pa u^j} \label{a1}
\ee
and
\be
\mu^{\alpha \beta}
\omega_{\alpha, ij} (u)
\omega_{\beta, ks} (u) =
\mu^{\alpha \beta}
\omega_{\alpha, ik} (u)
\omega_{\beta, js} (u).  \label{b1}
\ee
Moreover,
\bea
&&
{\pa^2 b_{\alpha} \over \pa u^i \pa u^l} =
\rho \, \eta^{kj} {\pa \omega_{\alpha, ij} \over \pa u^l}
 {\pa a \over \pa u^k} +
\rho \, \eta^{kj} \omega_{\alpha, ij} (u)
{\pa^2 a \over \pa u^k \pa u^l} = \nn\\
&&
= \rho \, \eta^{kj} {\pa \omega_{\alpha, ij} \over \pa u^l}
 {\pa a \over \pa u^k} +
\rho \, \eta^{kj} \omega_{\alpha, ij} (u)
 \lambda \, \mu^{\gamma \beta}
\omega_{\gamma, kl} (u)  b_{\beta} (u) = \nn\\
&&
= \rho \, \eta^{kj} {\pa \omega_{\alpha, lj} \over \pa u^i}
 {\pa a \over \pa u^k} +
\rho \, \eta^{kj} \omega_{\alpha, lj} (u)
 \lambda \, \mu^{\gamma \beta}
\omega_{\gamma, ki} (u)  b_{\beta} (u),
\eea
whence we have
\be
{\pa \omega_{\alpha, ij} \over \pa u^l} =
{\pa \omega_{\alpha, lj} \over \pa u^i} \label{a2}
\ee
and
\be
\eta^{kj} \omega_{\alpha, ij} (u)
\omega_{\gamma, kl} (u) = \eta^{kj} \omega_{\alpha, lj} (u)
\omega_{\gamma, ki} (u). \label{b2}
\ee
It follows from (\ref{a1}) and (\ref{a2}) that
there exist locally some functions $\psi_{\alpha} (u),$
$1 \leq \alpha \leq L,$ such that
\be
\omega_{\alpha, ij} (u) = {\pa^2 \psi_{\alpha} \over \pa u^i \pa u^j}
\ee
and then the relations (\ref{b1}) and (\ref{b2}) are
equivalent to the nonlinear system (\ref{g3}), (\ref{r3}) for the
functions $\psi_{\alpha} (u),$ $1 \leq \alpha \leq L.$

In arbitrary local coordinates, we obtain the following
integrable description of all $N$-dimensional flat torsionless submanifolds
in $(N + L)$-dimensional
pseudo-Euclidean spaces.

\bt  For each $N$-dimensional flat torsionless submanifold
in an $(N + L)$-dimensional
pseudo-Euclidean space with a flat first fundamental form $g_{ij} (u)$
there locally
exist functions $\psi_{\alpha} (u),$ $1 \leq \alpha \leq L,$ such that
the second fundamental forms have the form
\be
(w_{\alpha})_{ij} (u) = \nabla_i \nabla_j \psi_{\alpha}, \label{f1}
\ee
where $\nabla_i$ is the
covariant differentiation defined by the Levi-Civita connection
generated by the metric $g_{ij} (u)$.
The class of $N$-dimensional flat torsionless submanifolds
in $(N + L)$-dimensional
pseudo-Euclidean spaces is described by the following
integrable system of nonlinear equations for
functions $\psi_{\alpha} (u),$ $1 \leq \alpha \leq L:$
\be
 \sum_{n=1}^N  \nabla^n \nabla_i \psi_{\alpha}
\nabla_n \nabla_l \psi_{\beta}
=  \sum_{n=1}^N  \nabla^n \nabla_i \psi_{\beta}
\nabla_n \nabla_l \psi_{\alpha}, \label{non1}
\ee
\be
\sum_{\alpha = 1}^L \sum_{\beta = 1}^L \mu^{\alpha \beta}
\nabla_i \nabla_j \psi_{\alpha} \nabla_k \nabla_l \psi_{\beta} =
\sum_{\alpha = 1}^L \sum_{\beta = 1}^L \mu^{\alpha \beta}
\nabla_i \nabla_k
\psi_{\alpha} \nabla_j \nabla_l \psi_{\beta} , \label{non2}
\ee
where $\nabla_i$ is the
covariant differentiation defined by the Levi-Civita connection
generated by a flat metric
$g_{ij} (u)$, $\nabla^i = g^{is} (u) \nabla_s,$
$g^{is} (u) g_{sj} (u) = \delta^i_j.$
Moreover, in this case
the systems
of hydrodynamic type
\be
u^i_{t_{\alpha}} = \left (\nabla^i \nabla_j
\psi_{\alpha}\right ) u^j_x,
\ \ \ \ 1 \leq \alpha \leq L, \label{h1}
\ee
are commuting integrable bi-Hamiltonian systems
of hydrodynamic type.
\et

Now we will also find some natural and very important
integrable reductions of the
nonlinear system (\ref{g3}), (\ref{r3}).

\section{Reduction to the associativity equations \\ of
two-dimensional topological quantum field \\ theories
and potential flat torsionless\\ submanifolds
in pseudo-Euclidean spaces} \label{section3}

\bt
If we take $L = N$, $\mu^{ij} = c \eta^{ij},$ $1 \leq i, j \leq N,$
$c$ is an arbitrary nonzero constant,
and $\psi_{\alpha} (u) = {\pa \Phi / \pa u^{\alpha}},$
$1 \leq \alpha \leq N,$ where $\Phi = \Phi (u^1, \ldots, u^N)$,
then the Gauss equations {\rm (\ref{g3})} coincide with the Ricci equations
{\rm (\ref{r3})} and both of them coincide with the associativity
equations of two-dimensional topological quantum field theories
for the potential $\Phi (u)$:
\be
\sum_{i = 1}^N \sum_{j = 1}^N
\eta^{ij}
\left ({\pa^3 \Phi \over \pa u^i \pa u^m \pa u^k}
{\pa^3 \Phi \over \pa u^j \pa u^n \pa u^l} -
{\pa^3 \Phi \over \pa u^i \pa u^m \pa u^l}
{\pa^3 \Phi \over \pa u^j \pa u^n \pa u^k} \right ) = 0,  \label{ass2}
\ee
\et
\bt
The associativity
equations of two-dimensional topological quantum field theories
describe a special class of $N$-dimensional
flat submanifolds without torsion
in $2N$-dimensional
pseudo-Euclidean spaces (a class of {\rm potential flat
torsionless submanifolds}).
\et
\bd
An $N$-dimensional flat torsionless submanifold in
a $2N$-dimensional pseudo-Euclidean space with a flat
first fundamental form $g_{ij} (u) d u^i d u^j$ is called {\it potential}
if there exist a function $\Phi (u)$ on the submanifold such that
the second fundamental forms of the submanifold have the form
\be
(\omega_i)_{jk} (u) d u^j d u^k =
\left ( \nabla_i \nabla_j \nabla_k \Phi (u) \right ) d u^j d u^k,
\ \ \ \ 1 \leq i \leq N,
\ee
where $\nabla_i$ is the
covariant differentiation defined by the Levi-Civita connection
generated by the flat metric
$g_{ij} (u)$.
\ed

According to the Bonnet theorem
any solution $\Phi (u)$
of the associativity equations
(with the corresponding constant metric $\eta_{ij}$)
defines a unique (up to motions) $N$-dimensional potential flat torsionless
submanifold of
the corresponding $2N$-di\-men\-si\-o\-nal pseudo-Euclidean space with the
first fundamental form $\eta_{ij} d u^i d u^j$ and the
second fundamental forms
$\omega_n (u) = (\partial^3 \Phi / (\partial u^n \partial u^i \partial u^j))
d u^i d u^j$
 given by the third derivatives of the
potential $\Phi_{\alpha} (u)$.
We do not distinguish solutions of
the associativity equations up to
terms quadratic in the coordinates $u$.

\bt
On each potential flat torsionless submanifold
in a pseudo-Euclidean space there is a structure of Frobenius
algebra given by the Weingarten operators $(A_s)^i_j (u) = -
\eta^{ik} (\omega_s)_{kj} (u)$,
\be
c^k_{ij} (u^1, \ldots, u^N) = \eta^{ks}
(\omega_i)_{sj} (u^1, \ldots, u^N).
\ee
In arbitrary local coordinates, the Frobenius structure is given by
\be
c^k_{ij} (u^1, \ldots, u^N) = g^{ks} (u^1, \ldots, u^N)
(\omega_i)_{sj} (u^1, \ldots, u^N),
\ee
where $g^{ij} (u)$ is the contravariant metric inverse to the
first fundamental form $g_{ij} (u)$,
$g^{is} (u) g_{sj} (u) = \delta^i_j,$
$(\omega_k)_{ij} (u) d u^i d u^j,$ $1 \leq k \leq N,$ are
the second fundamental forms.
\et

\bt
Each $N$-dimensional Frobenius manifold
can locally be represented as a potential flat torsionless
$N$-dimensional
submanifold in a $2N$-dimensional
pseudo-Euclidean space. This submanifold is
uniquely determined up to motions.
\et



\end{document}